\documentclass[10pt]{article}

\usepackage{color}
\usepackage{graphics}
\usepackage{latexsym}
\usepackage{amsmath}
\usepackage{amssymb}
\usepackage{pstricks}
\usepackage{pst-node}

\makeatletter
\newtheorem{theorem}{Theorem}[section]
\newtheorem{lemma}[theorem]{Lemma}

\newtheorem{remark}[theorem]{Remark}
\newtheorem{cor}[theorem]{Corollary}
\newtheorem{prop}[theorem]{Proposition}
\newtheorem{defin}[theorem]{Definition}

\def\ps@headings{
 \def\@oddhead{\footnotesize\rm\hfill\runningheadodd\hfill\thepage}
 \def\@evenhead{\footnotesize\rm\thepage\hfill\runningheadeven\hfill}
 \def\@oddfoot{}
 \def\@evenfoot{\@oddfoot}
}

\newcommand{\Prf}{\noindent{\bf Proof}.\quad }

\newcommand{\qed}{\hfill$\Box$}

\makeatother

\def\runningheadeven{Adjacency Matrices of Configuration Graphs}
\def\runningheadodd{M. Abreu, M. Funk, D. Labbate, V. Napolitano}

\title{Adjacency Matrices of Configuration Graphs}
\author{{\rm M. Abreu,${}^\dagger$\thanks{This research was carried out within the activity
 of INdAM-GNSAGA and supported by the Italian Ministry MIUR}   $\;$
 M.J. Funk,$^{\dagger \ast }$$\;$D. Labbate,$^{\ddagger \ast}$ $\;$V. Napolitano$^{\S \ast }$}\\
\\ ${}^\dagger$ \small Dipartimento di Matematica e Informatica -
Universit\`a della
        Basilicata \\
      \small 85100 Potenza -
         Italy (marien.abreu@unibas.it; martin.funk@unibas.it)
          \\
\\ $^\ddagger$ \small Dipartimento di Matematica - Politecnico di Bari
\\ \small Via E. Orabona, 4 - 70125 Bari -
         Italy (labbate@poliba.it)\ \\
         \\ $^\S$ \small Dipartimento di Ingegneria Civile - Seconda Universit\`a di Napoli
\\ \small  via Roma, 29 - 81031 Aversa -
         Italy (vito.napolitano@unina2.it)
         }

      \date{}

\begin{document}
\maketitle
\pagestyle{headings}

\begin{abstract}
In 1960, Hoffman and Singleton \cite{HS60} solved a celebrated
equation for square matrices of order $n$, which can
be written as
$$ (\kappa - 1) \; I_n + J_n - A A^{\rm T}\;=\; A$$
where  $I_n$, $J_n$, and $A$ are the identity matrix, the all one
matrix, and a $(0,1)$--matrix with all row and column sums equal
to $\kappa$, respectively. If $A$ is an incidence matrix of some
configuration $\cal C$ of type $n_\kappa$, then the left-hand side
$\Theta(A):= (\kappa - 1)I_n + J_n - A A^{\rm T}$ is an
adjacency matrix of the non--collinearity graph $\Gamma$ of $\cal C$.
In certain situations, $\Theta(A)$ is also an incidence matrix of some $n_\kappa$ configuration,
namely the neighbourhood geometry of $\Gamma$ introduced by
Lef\`evre-Percsy, Percsy, and Leemans \cite{LPPL}.

The matrix operator $\Theta$ can be reiterated and we pose the problem of solving
the generalised Hoffman--Singleton equation $\Theta^m(A)=A$. In particular,
we classify all $(0,1)$--matrices $M$ with all row and column sums equal to $\kappa$, for $\kappa = 3,4$,
which are solutions of this equation.
As a by--product, we obtain characterisations for incidence matrices of the configuration $10_3F$
in Kantor's list \cite{Kantor} and the $17_4$ configuration $\#1971$ in Betten and Betten's list \cite{BB99}.
\end{abstract}

\section{Preliminaries: $(0,1)$--Matrices}\label{matr}

Denote by  $I_n$ and $J_n$ the identity matrix and the all one matrix of order $n$,
respectively. A $(0,1)$--matrix is said to be $J_2$--{\em free}  if it does not contain an all one submatrix
of order $2$. In order to clearly display $(0,1)$--matrices, we will often omit the entries $0$.

For integers $n$ and $\kappa$ ranging over and $\mathbb{N^+}$ and $\mathbb{Z}$ respectively, define
$$\delta(\kappa) := n - \kappa^2 + \kappa -1 \,.$$
Denote by ${\frak Z}_{n,\kappa}$ and  ${\frak  D}_{n,\kappa}$ the subclasses of all matrices in ${\frak M}_n(\mathbb{Z})$ and
${\frak M}_n(0,1)$, respectively,
which have all row and column sums equal to $\kappa$. Clearly, ${\frak
 D}_{n,\kappa}$ is empty for $\kappa < 0$. Put
${\frak Z}_{n} := \cup_{\kappa} {\frak Z}_{n,\kappa}$ and ${\frak
D}_{n} := \cup_{\kappa} {\frak D}_{n,\kappa}$. A key r\^ ole will
be played by the application
$$\Theta : {\frak Z}_{n}  \;\longrightarrow \; {\frak M}_n(\mathbb{Z}) \quad \makebox{defined by} \quad \Theta (A)
\;:=\;(\kappa - 1) \; I_n + J_n - A A^{\rm T} $$ where, as usual, $A^{\rm T}$ denotes the transpose of $A$.

    \begin{lemma}\label{matrixop}
The application $\Theta$ is a matrix operator in the class ${\frak Z}_{n}$.
    \end{lemma}

\Prf We show that $A=(a_{ij}) \in {\frak Z}_{n,\kappa}$  implies $\Theta(A)
\in {\frak Z}_{n,\delta(\kappa)}$.  Then the $i^{th}$ row sum and the $j^{th}$ column sum
of $A A^{\rm T}$ respectively read
$$\sum_{j = 1}^{n} \sum_{k = 1}^{n} a_{ik}a_{jk} = \sum_{k =
1}^{n} a_{ik}(\sum_{j = 1}^{n} a_{jk})=  \sum_{k = 1}^{n}  a_{ik}
\,\kappa= \kappa^2 \quad \makebox{and}$$$$\sum_{i = 1}^{n} \sum_{k
= 1}^{n} a_{ik}a_{jk} = \sum_{k = 1}^{n} (\sum_{i = 1}^{n}a_{ik})
a_{jk} = \sum_{k = 1}^{n} \kappa\, a_{jk} = \kappa^2\,.$$
Hence all row and column sums of $A A^{\rm T}$ have
constant  value $\kappa^2$. This, in turn, implies that the
summands $(\kappa-1) I_n$, $J_n$ , and $- A A^{\rm T}$ contribute
$\kappa -1$, $n$ and $-\kappa^2$, respectively, to each row and
column sum of $\Theta(A)$. \qed

\

For positive integers $m$ and $\kappa$, we are interested in the
subclass ${\frak S}_{\kappa} \subseteq {\frak D}_{\kappa^2+1,\kappa}$ of
solutions for the {\em generalised Hoffman--Singleton} matrix
equation
$$
(gHS) \qquad\qquad\qquad\qquad\qquad \qquad\Theta^m(A) \; = \;
A\,. \qquad \qquad \qquad \qquad\qquad \qquad
$$

The proof of Lemma \ref{matrixop} shows that $A \in {\frak D}_{n,\kappa}$
is a solution of $(gHS)$ only if $\delta^m(\kappa)
= \kappa$. This Diophantine equation in $n,m,\kappa$ has solutions $(\kappa^2+1,m, \kappa)$ for
$m \in \mathbb{N}^+$ and $\kappa \in \mathbb{Z}$
and in this paper we will consider this case only.
We conjecture that {\em the only other solutions to this equation are
$((\kappa-1)^2+2,2\mu, \kappa)$ for $\mu \in \mathbb{N}^+$ and $\kappa \in \mathbb{Z}$.}

\begin{lemma}\label{lezero}
Let $d_{ii}$ be an entry on the main diagonal of $\Theta(A)$ for
some $A \in {\frak Z}_{n,\kappa}$. Then $d_{ii} \le 0$ and
equality holds if and only if $A \in {\frak D}_{n,\kappa}$.
\end{lemma}

\Prf Each entry on the main diagonal of $A A^{\rm T}$, say
$b_{ii}$, is the product of the $i^{th}$ row of $A$ with itself.
Hence $b_{ii} = \sum_{j=1}^n a_{ij}^2$ is a sum of squares.
Collect the positive summands of  $\sum_{j=1}^n a_{ij} = \kappa$
and write it as $$\sum_{j=1}^n a_{ij} = \sum_{k \in K} p_k -
\sum_{l \in L} n_l$$ where $p_k \ge 1$  and $n_l \ge 0$ are non
negative integers for suitable (possibly empty) index sets $K$
and $L$, respectively. Since $p_k^2 \ge p_k$ and $n_l^2 \ge n_l$,
one has
$$b_{ii} = \sum_{j=1}^n a_{ij}^2 = \sum_{k \in K} p_k^2 + \sum_{l \in L}
n_l^2 \ge \sum_{k \in K} p_k -
 \sum_{l \in L} n_l = \sum_{j=1}^n a_{ij} = \kappa\,.$$
 Thus the summands $(\kappa-1) I_n$, $J_n$, and $- A A^{\rm T}$
contribute $\kappa -1$, $1$ and $-b_{ii}$, respectively, to each
entry $d_{ii}$ on the main diagonal of $\Theta(M)$. This implies
 $d_{ii} = \kappa - 1 + 1 - b_{ii} = \kappa - b_{ii} \le
0$.

Equality holds if and only if
 $$\sum_{k \in K} p_k^2 -  \sum_{k \in K} p_k= - \sum_{l \in L} n_l^2
-  \sum_{l \in L} n_l \,.$$ Since the left-hand side $\sum_{k \in
K} (p_k^2 -  p_k)$ is either zero or positive, whereas the
right-hand side $ - \sum_{l \in L} (n_l^2 + n_l)$ is either zero
or negative, both sides must be zero. Clearly, the right-hand side
is zero if and only if $n_l = 0$ for all $l \in L$, whereas the
left-hand side is zero if and only if $p_k^2 = p_k$ and thus
$p_k=1$ for all $k \in K$.  Hence $b_{ii}$ assumes its minimum value,
namely $\kappa$, if and only if the $i^{th}$ row of $A$ is made up
of $\kappa$ entries $1$ and $n-\kappa$ entries $0$.
\qed
\medskip

\begin{remark}
By Lemma \ref{matrixop}, the matrix $\Theta^m(A)$ has
all row and column sums equal to $\delta^m(\kappa)$. Hence
$\delta^m(\kappa) = \kappa$ for $A \in {\frak S}_{\kappa}$.
For $\kappa = 1$ the two solutions of the Diophantine equation $\delta^m(\kappa) = \kappa$ coincide.
The second solutions $((\kappa-1)^2+2,2\mu, \kappa)$ %for $\mu \in \mathbb{N}^+$ and $\kappa \in \mathbb{Z}$,
do not contribute to our problem for $\kappa \ge 2$.
In fact, for $n = (\kappa-1)^2 + 2$,
the value of $\delta(\kappa) = 2 - \kappa$ is no longer positive.
Appling Lemma \ref{lezero} to
$\Theta(\Theta(A)) = \Theta^2(A)$, each entry on the main diagonal
of $\Theta^2(A)$ is negative, thus $\Theta^2(A) \ne A$.
Hence, by induction we get the contradiction $\Theta^{2\mu}(A) \ne A$
for all $\mu \in \mathbb{N}^+$.
\end{remark}

\begin{theorem}\label{dlsz}
Let $A \in {\frak S}_{\kappa}$ be a solution of $\Theta^m(A)=A$, $m\geq 1$. Then:
$$
\begin{array}{rcl}
(L)&& A \; \makebox{is}\; J_2\makebox{--free}\,.\\
%(O) &&  A \; \makebox{has order  }n = \kappa^2 + 1\,;\\
(S) && A \; \makebox{is symmetric;}\\
(Z) && A \; \makebox{has entries } 0  \makebox{ on its main
diagonal;}\\
\end{array}$$
\end{theorem}

\Prf  $(L)$ By hypothesis $\Theta^m(A) = A$, for some $m \ge 1$.
If $A$ contained entries $1$ in positions $(i,j)$, $(i,k)$,
$(l,j)$, and $(l,k)$ for $i,j,k,l \in \{1, \ldots, n\}$ with $i
\ne l$ and $j \ne k$, then the entry $(A A^{\rm T})_{il} = \sum_{r
= 1}^n a_{ir}a_{lr}$ would have at least two summands $1$, namely
for $r = j$ and $r = k$. This would imply $(A A^{\rm T})_{il} \ge
2$ and $(\Theta(A))_{il} \le - 1$, hence $\Theta(A) \ne A$. The
only remaining option would be $\Theta^m(A) = A$ for some $m \ge
2$. On the other hand, $(\Theta(A))_{il} \le - 1$ means $\Theta(A)
\not \in {\frak D}_{\kappa^2+1,\kappa}$. Then Lemma \ref{lezero}
implies that each entry on the main diagonal of $\Theta^2(A)$ is
negative, hence $\Theta^2(A) \not \in {\frak D}_{\kappa^2+1,\kappa}$.
Induction on $m$ shows $\Theta^m(A) \not \in {\frak D}_{\kappa^2+1,\kappa}$
for all $m \ge 2$, a contradiction.

$(S)$ Defined as a sum of symmetric matrices, $\Theta(A)$ is
symmetric for any $A \in {\frak M}_n(\mathbb{Z})$.
Hence $A \in {\frak S}_{\kappa}$ implies that $A =
\Theta(\Theta^{m-1}(A))$ is symmetric.

$(Z)$ follows immediately from Lemma \ref{lezero}.
 \qed
\medskip

Note that conditions $(L)$, %$(O)$,
$(S)$, and $(Z)$ do not
characterise the class $\frak{S}_\kappa$.
A counterexample will be presented in Remark
\ref{conterex}$(ii)$.

\section{Connection to Configurations and Graphs}

For notions from graph theory and incidence geometry, we
respectively refer to \cite{BondyMurty} and \cite{Dembowski}. We
consider undirected graphs without loops or multiple edges.

A graph is said to be $C_4$--{\em free} if it does not contain $4$--cycles.
With each permutation $\pi$ in the symmetric group ${\cal
S}_n$ we can associate its {\it permutation matrix} $P_\pi =
(p_{ij})_{1 \le i,j \le n}$ which is defined by $p_{ij} = 1$ if
$i^{\,\pi} = j$, and $p_{ij} = 0$ otherwise.

We call an incidence structure (in the sense of \cite{Dembowski})
{\em linear} if any two distinct points are incident with at most
one line. A {\em  configuration $\cal C$ of type} $n_\kappa$ is a
linear incidence structure consisting of $n$ points and $n$ lines
such that each point and line is  incident with $\kappa$ lines and
points, respectively. To individualise certain $10_3$ and
$17_4$ configurations, we refer to the lists in \cite{Kantor}
and \cite{BB99}, respectively.

Fix a labelling for the points and lines of a configuration $\cal
C$ and consider the {\it incidence matrix} $C$ of $\cal C$ (cf.
e.g. \cite[pp. 17--20]{Dembowski}): there is an entry $1$ and $0$
in position $(i,j)$ of $C$ if and only if the point $p_i$ and the
line $l_j$ are incident and non--incident, respectively. The
following result is well known.

    \begin{lemma}\label{confprop}

$(i)$ A square $(0,1)$--matrix $C$ of order $n$ is an
incidence matrix of some configuration $\cal C$ of type $n_\kappa$
if and only if it is $J_2$--free and has all row and column sums
equal to $\kappa$.

$(ii)$ Any other incidence matrix of $\cal C$ has the form
$S_l^{-1}CS_r$ for permutation matrices $S_l$ and $S_r$ of order
$n$, corresponding to re--labellinga of points and lines.

$(iii)$ $C$ is symmetric and has entries $0$ on its main diagonal if and only if
$\,\cal C$ admits a self--polarity $p_i \longleftrightarrow l_i$
without absolute elements.\qed
    \end{lemma}

Adjacency matrices depend on  the labelling of the vertices. If
$G$ is a graph of order $n$ and $A$ an adjacency matrix for $G$,
then any other adjacency matrix has the form $S^{-1}AS$ for a
suitable permutation matrix $S$ of order $n$ (which represents a
re--labelling of the vertices).

Any two matrices $M_1$ and $M_2$
of order $r$ are said to be {\em permutationally equivalent} or {\em
p--equivalent} for short, denoted by $M_1 \sim M_2$, if there exists
a permutation matrix $S$ of order $r$ such that $M_2 =
S^{-1}M_1S$. The following result is again well known.

\begin{lemma}\label{adj}
$(i)$ A square $(0,1)$--matrix $A$ of order $n$ is an adjacency matrix
of some $\kappa$--regular graph $G$ of order $n$ if and only if it
is symmetric, has entries $0$ on its main diagonal, and all row
and column sums equal to $\kappa$.

$(ii)$ Any other adjacency matrix of $G$ is
p--equivalent to $A$.

$(iii)$ A graph $G$ is $C_4$--free if and only if its
adjacency matrix is $J_2$--free. \qed
\end{lemma}

\

With each configuration $\cal C$ of type $n_\kappa$, we can
associate its {\it configuration graph} $\Gamma({\cal C})$, known
also as {\it non--collinearity graph}, as the result of the following operation $\Gamma$:
the vertices of $\Gamma({\cal C})$ are the points of $\cal C$; any two vertices
are joined by an edge if they are not incident with one and the
same configuration--line (\cite{Gropp93}). The number of points in
$\cal C$ not joined with an arbitrary point of $\cal C$ is given
by $\delta(\kappa) : = n - \kappa^2 + \kappa - 1$, called the {\it
deficiency} of $\cal C$. Finite projective planes are
characterised by deficiency $0$. Thus the configuration graph
$\Gamma({\cal C})$ is a $\delta(\kappa)$--regular graph on $n$
vertices. Since the following Lemma plays a key r\^ ole, we will
also quote its short proof: %from \cite{AFLN1}.

  \begin{lemma}\label{formula}\cite{AFLN1}
Let $\mathcal C$ be a configuration of type $n_\kappa$ with
incidence matrix $C$. Then the adjacency matrix $A$ of the
configuration graph $\Gamma({\cal C})$ is given by
$$
A \, = \, (\kappa - 1) \; I_n + J_n - C \;C^{\rm T} \,.
$$
    \end{lemma}

\Prf Let $M :=(m_{i,j}) := C \;C^{\rm T}$. An arbitrary entry
$m_{i,j}$ of $M$ is the result of the usual dot product (over
$\mathbb{R}$) of the $i^{th}$ row and the $j^{th}$ row of $C$.
Since the rows represent the $i^{th}$ and $j^{th}$ points of $\cal
C$, say $p_i$ and $p_j$, we have
$$ m_{i,j} \quad =
\quad \left\{\begin{array}{cl} \kappa & \makebox{if  } \; i = j
\,;
\\ 1 & \makebox{if } \; i \ne j\;\; \makebox{ and there is a
line in} \;\,{\cal C}\,\;\makebox{joining} \;\; p_i , p_j\,;\\
0 & \makebox{if } \; i \ne j\; \;\makebox{ and}\;\; p_i,
p_j\;\;\makebox{are not joined by any line of} \;\,{\cal C}.
\end{array}\right.$$
On the other hand, the adjacency matrix $A := (a_{i,j})$ of the
configuration graph $\Gamma({\mathcal C})$ has entries:
$$
a_{i,j} \quad = \quad \left\{\begin{array}{cl} 1 & \makebox{if the
points}\;\; p_i , p_j \;\; \makebox{are not joined by any line of}
\;\,{\cal C}\,;
\\ 0 & \makebox{otherwise}\;.
\end{array}\right.$$
This implies $A \;=\;(\kappa - 1) \;  I_n + J_n - M $. \qed
\medskip

Recently, Lef\`evre-Percsy,  Percsy, and Leemans \cite{LPPL}
introduced an operation $\cal N$ which can be seen as a kind of
\lq \lq inverse" operation for $\Gamma$. It associates with each
graph $G$ its {\em neighbourhood geometry} ${\cal N}(G) =
(P,B,|)$: let $P$ and $B$ be two copies of $V(G)$, whose
elements are called {\it points} and {\it blocks}, respectively; a
point $x \in P$ is incident with a block $b \in B$ (in symbols $x
| b)$ if and only if $x$ and $b$, seen as vertices in $G$, are
adjacent. On the other hand, $\cal N$ has no effect in terms of $(0,1)$--matrices.
In fact, it only reinterprets an adjacency matrix of $G$ as an incidence
matrix of some configuration, namely ${\cal N}(G)$.

\begin{remark}\label{conterex}

$(i)$ Given a configuration $\cal C$ of type $(\kappa^2 +1)_\kappa \, $, its
configuration graph $\Gamma({\cal C})$ need not be $C_4$--free
(i.e. it may contain a $4$--cycle $p_1,l_1,p_2,l_2$). If this happens,
the neighbourhood geometry of $\Gamma({\cal C})$ contains a di-gon
$(\,\{p_1,p_2\}\,,\, $ $\{l_1,l_2\}\,, \, \{(p_i,l_j) \;|\,i,j =
1,2\})$ which is forbidden for a linear incidence structure.
Hence we say that {\em $\Gamma({\cal C})$ is not $\cal N$--admissible}.

$(ii)$ A $C_4$--free $\kappa$--regular graph $G$ on $\kappa^2 +1$ vertices
is ($\Gamma \circ {\cal N}$)--{\em admissible}, whereas a $(\kappa^2 +1)_\kappa$
configuration $\cal C$ is said to be (${\cal N} \circ \Gamma$)--{\em admissible}
if its configuration graph $\Gamma (\cal C)$ is $C_4$--free (cf. $(i)$).

$(iii)$ Recall that a {\em Terwilliger graph} is a non--complete graph $G$ such that,
for any two vertices $v_1, v_2 \in V(G)$ at distance $2$ from each other,
the induced subgraph $G[N_G(v_1) \cap N_G(v_2)]$ is a clique of size $\mu$, for some fixed $\mu \ge 0$
(cf. e.g. {\rm \cite[p.  34]{BCN}}).
Thus the class of $(\Gamma \circ {\cal N})$--admissible $\kappa$--regular graphs coincides
with the class of $\kappa$--regular Terwilliger graphs for $\mu =1$.

$(iv)$ The conditions $(L)$,
%$(O)$,
$(S)$, and $(Z)$ of Theorem \ref{dlsz} do not characterise the class $\frak{S}_\kappa$.
A counterexample is given by the following adjacency matrix $A_1$ of the Terwilliger graph $T_1$ since $A_1 \ne
\Theta(A_1)$, but $\Theta(A_1) = \Theta^m(A_1)$ for all $m \ge 1$.
$$
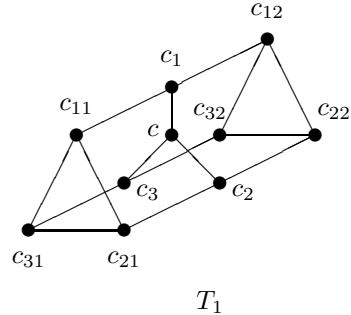
\begin{picture}(20,7)

\put(1.8,0){\makebox(0,0){$c_{32}$}}
\put(1.8,0.6){\makebox(0,0){$c_{22}$}}
\put(1.8,1.2){\makebox(0,0){$c_{12}$}}
\put(1.8,1.8){\makebox(0,0){$c_{31}$}}
\put(1.8,2.4){\makebox(0,0){$c_{21}$}}
\put(1.8,3){\makebox(0,0){$c_{11}$}}
\put(1.8,3.6){\makebox(0,0){$c_{3}$}}
\put(1.8,4.2){\makebox(0,0){$c_{2}$}}
\put(1.8,4.8){\makebox(0,0){$c_{1}$}}
\put(1.8,5.4){\makebox(0,0){$c\;\,$}}
\put(9.25,6){\makebox(0,0){$c_{32}$}}
\put(8.5,6){\makebox(0,0){$c_{22}$}}
\put(7.75,6){\makebox(0,0){$c_{12}$}}
\put(7,6){\makebox(0,0){$c_{31}$}}
\put(6.25,6){\makebox(0,0){$c_{21}$}}
\put(5.5,6){\makebox(0,0){$c_{11}$}}
\put(4.75,6){\makebox(0,0){$c_{3}$}}
\put(4,6){\makebox(0,0){$c_{2}$}}
\put(3.25,6){\makebox(0,0){$c_{1}$}}
\put(2.5,6){\makebox(0,0){$c$}} \put(2.25,-.25){\line(0,1){6.5}}
\put(9.7,-.25){\line(0,1){6.5}} \put(2.8,-.25){\line(0,1){6.5}}
\put(5.05,-.25){\line(0,1){6.5}} \put(7.35,-.25){\line(0,1){6.5}}
\put(1.4,5.7){\line(1,0){8.3}} \put(1.4,-.25){\line(1,0){8.3}}
\put(1.4,1.55){\line(1,0){8.3}} \put(1.4,3.35){\line(1,0){8.3}}
\put(1.4,5.15){\line(1,0){8.3}} \put(3.25,5.4){\makebox(0,0){$
1$}} \put(4,5.4){\makebox(0,0){$
1$}}\put(4.75,5.4){\makebox(0,0){$ 1$}}
\put(2.5,4.8){\makebox(0,0){$ 1$}} \put(2.5,4.2){\makebox(0,0){$
1$}}\put(2.5,3.6){\makebox(0,0){$ 1$}}
\put(5.5,4.8){\makebox(0,0){$ 1$}} \put(6.25,4.2){\makebox(0,0){$
1$}}\put(7,3.6){\makebox(0,0){$ 1$}}
\put(7.75,4.8){\makebox(0,0){$ 1$}} \put(8.5,4.2){\makebox(0,0){$
1$}}\put(9.25,3.6){\makebox(0,0){$ 1$}}
\put(3.25,3){\makebox(0,0){$ 1$}} \put(4,2.4){\makebox(0,0){$
1$}}\put(4.75,1.8){\makebox(0,0){$ 1$}}
\put(3.25,1.2){\makebox(0,0){$ 1$}} \put(4,.6){\makebox(0,0){$
1$}}\put(4.75,0){\makebox(0,0){$ 1$}}

\put(6.25,3){\makebox(0,0){$ 1$}} \put(7,3){\makebox(0,0){$
1$}}\put(5.5,2.4){\makebox(0,0){$ 1$}} \put(7,2.4){\makebox(0,0){$
1$}} \put(5.5,1.8){\makebox(0,0){$
1$}}\put(6.25,1.8){\makebox(0,0){$ 1$}}

\put(8.5,1.2){\makebox(0,0){$ 1$}} \put(9.25,1.2){\makebox(0,0){$
1$}}\put(7.75,.6){\makebox(0,0){$ 1$}}
\put(9.25,.6){\makebox(0,0){$ 1$}} \put(7.75,0){\makebox(0,0){$
1$}}\put(8.5,0){\makebox(0,0){$ 1$}}
\put(.5,3){\makebox(0,0){$A_1\,${\rm :}}}

\put(12,1.5){\circle*{.3}} \put(14,1.5){\circle*{.3}}
\put(14,2.5){\circle*{.3}} \put(16,2.5){\circle*{.3}}
\put(13,3.5){\circle*{.3}} \put(15,3.5){\circle*{.3}}
\put(16,3.5){\circle*{.3}} \put(18,3.5){\circle*{.3}}
\put(15,4.5){\circle*{.3}}\put(17,5.5){\circle*{.3}}

\put(12,1.5){\line(2,1){4}} \put(14,1.5){\line(2,1){4}}
\put(13,3.5){\line(2,1){4}}

\put(12,1.5){\line(1,2){1}} \put(14,1.5){\line(-1,2){1}}
\put(12,1.5){\line(1,0){2}}

\put(16,3.5){\line(1,2){1}} \put(18,3.5){\line(-1,2){1}}
\put(16,3.5){\line(11,0){2}}

\put(15,4.5){\line(0,-1){1}} \put(14,2.5){\line(1,1){1}}
\put(16,2.5){\line(-1,1){1}}

\put(17,6.1){\makebox(0,0){$c_{12}$}}
\put(15.8,4){\makebox(0,0){$c_{32}$}}
\put(18.4,4){\makebox(0,0){$c_{22}$}}
\put(15,5.1){\makebox(0,0){$c_1$}}
\put(14.6,3.6){\makebox(0,0){$c$}}
\put(13,4.1){\makebox(0,0){$c_{11}$}}
\put(12,.9){\makebox(0,0){$c_{31}$}}
\put(14,.9){\makebox(0,0){$c_{21}$}}
\put(16.5,2.3){\makebox(0,0){$c_2$}}
\put(14.5,2.3){\makebox(0,0){$c_3$}}
\put(15.8,0){\makebox(0,0){$T_1$}}
\end{picture}
$$
\medskip

Note that $\Theta(A_1)$ is an adjacency matrix for the Petersen
graph, see {\rm \cite[Proposition 4.2]{AFLN1}}.
\end{remark}

\section{Standard Forms}

Motivated by the next result, this section is dedicated to finding
some standard representatives within each class of p--equivalentity in $\frak S_{\kappa}$.

\begin{lemma}\label{sol}
Let $A$ be a solution of $(gHS)$. Then any p--equivalent matrix
$B$ is also a solution of $(gHS)$.
\end{lemma}

\Prf Suppose $B = S^{-1}AS$ for some permutation matrix $S$. Then,
in general,
$$\Theta(S^{-1}AS) =(\kappa-1)\,I_n + J_n - (S^{-1}AS)^2=
S^{-1}((\kappa-1)\,I_n + J_n - A^2)S =S^{-1}\Theta(A)S $$ and
hence
$$\Theta^m(S^{-1}AS) = \Theta(\Theta^{m-1}(S^{-1}AS)) =
\Theta(S^{-1}\Theta^{m-1}(A)S)=$$$$\qquad\qquad\qquad\qquad\qquad
\qquad\qquad=
S^{-1}\Theta(\Theta^{m-1}(A))S)=S^{-1}\Theta^m(A)S$$ by induction
on $m \ge 1$ and, in particular, $\Theta^m(B) = B$ if $\Theta^m(A)
= A$ for some $m \ge 1$. \qed

\

In the sequel we will use the following result %(cf.\cite[Proposition 2.7]{AFLN1})
which holds for all $\kappa \ge 2$:

    \begin{prop}\cite[Proposition 2.7]{AFLN1}
Let $G$ be $\kappa$--regular graph on $\kappa^2+1$ vertices whose
adjacency matrix fulfils condition $(L)$. Then $G$ has diameter $diam(G) \le 3$.
In particular, $diam(G) = 2$ if and only if $G$ has girth $5$.
\qed
    \end{prop}

Recall that a vertex $v$ of a graph $G$ is said to be a {\em
centre of $G$ with radius $2$} if the distance $d_G(v,w) \leq 2$
for each $w \in V(G)$.  In general, a graph $G$ with $diam(G)=3$
need not admit a centre with radius $2$, but we can prove the
following result, which also shows that the Conjecture
posed in \cite[p. 119]{AFLN1} holds true.

\begin{prop}\label{174}
Let $G$ be a $\kappa$--regular graph on $\kappa^2+1$ vertices
where $\kappa = 2,3,4$. If $G$ admits an adjacency matrix $A$
which is a solution of $(gHS)$, then $G$ has a centre with radius $2$.
\end{prop}

\Prf First we verify the following claim: {\it a vertex $v \in
V(G)$ does not lie in a $3$--cycle of $G$ if and only if $v $ is a
centre of $G$ with radius $2$.} To see this, let $v_1, \ldots,
v_\kappa$ denote the $\kappa$ vertices at distance $1$ from $v$.
Then $v$ does not lie in a $3$--cycle of $G$ if and only if we
encounter further $\kappa-1$ vertices $v_{ij}$, $j = 1, \ldots,
\kappa-1$, at distance $1$ from each $v_i$. (Since $A$ is
$J_2$--free by Theorem \ref{dlsz} and hence $G$ is $C_4$--free,
the vertices $v_{ij}$ turn out to be distinct in pairs.) This, in
turn, holds true if and only if there are no vertices at distance
$3$ from $v$ since $\{v,v_i,v_{ij} \,|\, i = 1, \ldots \kappa, j =
1, \ldots, \kappa-1\}$ is all of $V(G)$.

Secondly, a short calculation verifies that  $\kappa^2+1 \not
\equiv 0$ (mod $3$) for every integer $\kappa$. Hence $V(G)$
cannot be partitioned into vertex--disjoint $3$--cycles.

If $\kappa = 2$, then $G$ is a $5$-cycle and every vertex is a
centre with radius $2$. If $\kappa = 3$ and hence $|V(G)| = 10$,
then $G$ contains at most three disjoint $3$--cycles in $G$
and the remaining vertex is a centre with radius $2$.

Now let $\kappa =4$ and suppose that $G$ has no centre with radius
$2$. Then the above claim implies that $G$ contains at least one
vertex $v_0$ lying in two different $3$--cycles of $G$, say
$v_0v_1v_2v_0$ and $v_0v_3v_4v_0$. With respect to a labelling of
the vertices which starts with $v_0, v_1, v_2, v_3, v_4, \ldots$,
the first five rows and columns of the corresponding adjacency
matrix $B$ and $\Theta(B)$ read
$$\array{c|cccccc}
&v_0&v_1&v_2&v_3&v_4&\ldots\\\hline v_0&0& 1&1&1&1&\ldots\\
v_1&1&0&1&0&0&\ldots\\
v_2&1&1&0&0&0&\ldots\\v_3&1&0&0&0&1&\ldots\\v_4&1&0&0&1&0&\ldots\\
 \vdots&\vdots&\vdots&\vdots&\vdots&\vdots&\ddots\\
\endarray \quad\makebox{and}
\quad \left(\array{cccccc}
0&0&0&0&0&\ldots\\0&0&0&0&0&\ldots\\0&0&0&0&0&\ldots\\
0&0&0&0&0&\ldots\\0&0&0&0&0&\ldots\\
\vdots&\vdots&\vdots&\vdots&\vdots&\ddots\
\endarray \right) \,, $$ respectively. Denote by $S$ the submatrix of order $5
\times 12$ of $\Theta(B)$ made up by the first five rows and the
$6^{th}$ through the $17^{th}$ columns. Since $B$ is p--equivalent
to $A$ and $A$ is a solution of $(gHS)$, the image $\Theta(B)$ is
again an adjacency matrix of a $4$--regular graph on $17$
vertices (cf. Theorem \ref{dlsz}).  Thus $S$ must contain four
entries $1$ in each row. A short argument shows that this is not
compatible with Property $(L)$ of Theorem \ref{dlsz}, namely that
$\Theta(B)$ and hence $S$ is $J_2$--free. In fact, in the best
case we would need $13$ columns to fit four entries $1$ into each
row without producing a submatrix of $S$  isomorphic to $J_2$,
e.g.
$$\left(\array{cccccccccccccc}
1&1&1&1&&&&&&&&&\\1&&&&1&1&1&&&&&&\\
1&&&&&&&1&1&1&&&\\1&&&&&&&&&&1&1&1\\
&1&&&1&&&1&&&1&&\\
\endarray \right) \,. $$ This contradiction shows that $G$ has a
centre with radius $2$.
 \qed
\medskip

\begin{cor}\label{label} Let $G$ be a $\kappa$--regular graph on $\kappa^2+1$ vertices
where $\kappa = 2,3,4$ and let $v_0$ be a centre with radius $2$.
If $G$ admits an adjacency matrix $A$ which is a solution of
$(gHS)$, then the vertex set of $G$ is $$\{v_0\,,\;v_1\,,\ldots,
v_\kappa \,, \;v_{1,1}\,,\ldots,v_{1,\kappa-1}\;,\; v_{2,1}\,,
\ldots,v_{2,\kappa-1}\;, \;\ldots \;,\; v_{\kappa,1}\,, \ldots,
v_{\kappa,\kappa-1}\}$$  where,  $v_i$ and $v_{ij}$
denote the $\kappa$ neighbours
 of $v_0$ and  the $\kappa-1$ neighbours  of each
$v_i$ other than $v_0$, respectively,  for $i= 1,\ldots,\kappa$
and $j = 1, \ldots,\kappa-1$. \qed
\end{cor}

Denote by ${\bf 0}_\nu$  and ${\bf 1}_\nu$   row vectors of
dimension $\nu$ all of whose entries are $0$ and $1$,
respectively. Let ${\bf 0}_{\kappa,\kappa}$ be a copy of the zero
matrix of order $\kappa$ and $K = I_\kappa \otimes {\bf
1}_{\kappa-1}$ the Kronecker product of matrices. With respect to
the labelling mentioned in Corollary \ref{label}, the adjacency
matrix of $G$   gets the {\it standard form}
$$S(P)\quad := \quad \left(\begin{array}{llll}0 & {\bf 1}_\kappa
&& {\bf 0}_{\kappa^2-\kappa}\\
\rule{0pt}{14pt} {\bf 1}_\kappa^{\rm T}& {\bf 0}_{\kappa,\kappa} && K\\
\rule{0pt}{14pt} {\bf 0}_{\kappa^2-\kappa}^{\rm T} & K^{\rm T} && P
\end{array}\right)$$ where $P$ is a symmetric $(0,1)$--matrix of order
$\kappa^2-\kappa$ and all row and column sums equal to $\kappa-1$. We
regard $P$ as a block matrix $P = (P_{ij})_{1 \le i ,j,\le
\kappa}$ of order $\kappa$ with square blocks $P_{ij}$ of order
$\kappa-1$.

\begin{theorem} Let $A$ be a solution of $(gHS)$ and suppose
$\kappa \le 4$.

\noindent $(i)$ $A$ is p--equivalent to a standard form
$S(P)$ where each block $P_{ij}$ is a $(0,1)$--matrix which has at
most one entry $1$ in each row and column.

\noindent $(ii)$ If there exists a zero block in some row and
column of $P$, then all other blocks of $P$ in that row and column
are permutation matrices.

\noindent $(iii)$ If there exists a zero block on the diagonal of
$P$, then we can write $S(P)$ in such a way that $P_{11} = {\bf
0}_{\kappa-1,\kappa-1}$ and  $P_{1,i}= P_{i,1}=I_{\kappa-1}$ for
all $i=2,\ldots,\kappa$.
\end{theorem}

\Prf $(i)$ For $\kappa \le 4$, Corollary \ref{label} guarantees that $A$
has a p--equivalent standard form $S(P)$,  which is a solution of $(gHS)$.
Suppose that there are two entries $1$ in one and the same row
of some block $P_{i,j}$, say in positions $(s,t)$ and $(s,u)$, for some
$s,t,u \in \{1, \ldots, \kappa -1\}$. Then in $S(P)$, we encounter entries $1$ in
the following four positions:
$$
\begin{array}{l}
(1+j \, , \,  1+\kappa+(j-1)(\kappa-1)+t)\\
(1+j \, , \,  1+\kappa+(j-1)(\kappa-1)+u)\\
(1+\kappa+(i-1)(\kappa-1)+s  \, , \,   1+\kappa+(j-1)(\kappa-1)+t)\\
(1+\kappa+(i-1)(\kappa-1)+s  \, , \,  1+\kappa+(j-1)(\kappa-1)+u),
\end{array}
$$
thus obtaining a $J_2$ submatrix of $S(P)$, which contradicts Property $(L)$ of Theorem \ref{dlsz}.
By symmetry, an analogous reasoning works if there
were two entries $1$ in one and the same column of some block
$P_{ij}$.

$(ii)$ follows for arithmetic reasons:  $P$ is at once a
block matrix of order $\kappa$  and a $(0,1)$--matrix with
all row and column sums equal to $\kappa-1$.

$(iii)$ is the result of a suitable relabelling of the
rows and columns of $S(P)$. Consider the vertices $v_0, v_i, v_{ij}$
introduced in Corollary \ref{label}: if $P_{ii} = {\bf
0}_{\kappa-1,\kappa-1}$, then first exchange the r\^ oles of $v_1$
and $v_i$; secondly, for $i = 2, \ldots, \kappa$, relabel the
vertices within each family $\{v_{ij}\,|\, j = 1, \ldots,
\kappa-1\}$ such that $v_{1,j}$ is adjacent with $v_{ij}$ for all
$j = 1, \ldots, \kappa$. \qed

    \begin{defin}\label{HS}
The standard form $S(P)$ is said to be a {\em Hoffman--Singleton form},
or an {\em $HS$--form} for short, denoted by $S_{HS}(P)$, if all
the diagonal blocks $P_{ii}$ are zero blocks and
$P_{1,i}=P_{i,1}=I_{\kappa-1}$ for all $i=2,\ldots,\kappa$ {\rm
\cite{HS60}}.
    \end{defin}

\section{Classification}\label{class}

In 1960, Hoffman and Singleton \cite{HS60} classified all
$\kappa$--regular graphs on $\kappa^2+1$ vertices with girth $5$.
By eigenvalue techniques, they actually proved the following

    \begin{theorem}\label{resultfork} {\rm \cite{HS60}}
Let  $A \in {\frak D}_{\kappa^2+1,\kappa}$ be a solution of the
Hoffman--Singleton equation $\Theta(A) = A$. Then one of the following statements holds:
$$
\begin{array}{rcl}
(i) & \kappa = 2 & \makebox{and}\; A \; \makebox{is an adjacency
matrix for the 5--cycle};\\
(ii) & \kappa = 3 & \makebox{and}\; A \; \makebox{is an adjacency
matrix for the Petersen graph as well}\\
 &  & \makebox{as an incidence matrix for the Desargues configuration};\\
(iii) & \kappa = 7 & \makebox{and}\; A \; \makebox{is an adjacency
matrix for Hoffman--Singleton's graph}; \\(iv) & \kappa = 57 &
\makebox{(no graph or configuration is known). \qquad \qquad \qquad \qquad \qquad \qquad \qed}
\end{array}
$$
\end{theorem}

An  HS--form for $\kappa = 7$, i.e. for an adjacency matrix of the
Hoffman--Singleton graph, is presented in \cite[Figure 3]{HS60}.
The following two matrices are HS--forms for adjacency matrices of
the $5$--cycle and the Petersen graph, respectively:

    $$\left(\begin{array}{c|cc|c|c}
    &1&1&&\\\hline1&&&1&\\1&&&&1\\\hline&1&&&1\\\hline&&1&1&\\
    \end{array}\right) \qquad\qquad
    \left(\begin{array}{c|ccc|cc|cc|cc}
    &1&1&1&&&&&&\\\hline
    1&&&&1&1&&&&\\
    1&&&&&&1&1&&\\
    1&&&&&&&&1&1\\\hline
    &1&&&&&1&&1&\\
    &1&&&&&&1&&1\\\hline
    &&1&&1&&&&&1\\
    &&1&&&1&&&1&\\\hline
    &&&1&1&&&1&&\\
    &&&1&&1&1&&&\\
\end{array}\right)
    $$

In \cite{AFLN1}, two further solutions for the generalised
Hoffman--Singleton equation $\Theta^m(A) = A$, with $m \ge 2$, have been found;
the first being an adjacency matrix of the following Terwilliger graph $T_2$:

\begin{picture}(20,8)

\put(1.8,0){\makebox(0,0){$c_{32}$}}
\put(1.8,0.6){\makebox(0,0){$c_{31}$}}
\put(1.8,1.2){\makebox(0,0){$c_{22}$}}
\put(1.8,1.8){\makebox(0,0){$c_{21}$}}
\put(1.8,2.4){\makebox(0,0){$c_{12}$}}
\put(1.8,3){\makebox(0,0){$c_{11}$}}
\put(1.8,3.6){\makebox(0,0){$c_{3}$}}
\put(1.8,4.2){\makebox(0,0){$c_{2}$}}
\put(1.8,4.8){\makebox(0,0){$c_{1}$}}
\put(1.8,5.4){\makebox(0,0){$c\;\,$}}
\put(9.25,6){\makebox(0,0){$c_{32}$}}
\put(8.5,6){\makebox(0,0){$c_{31}$}}
\put(7.75,6){\makebox(0,0){$c_{22}$}}
\put(7,6){\makebox(0,0){$c_{21}$}}
\put(6.25,6){\makebox(0,0){$c_{12}$}}
\put(5.5,6){\makebox(0,0){$c_{11}$}}
\put(4.75,6){\makebox(0,0){$c_{3}$}}
\put(4,6){\makebox(0,0){$c_{2}$}}
\put(3.25,6){\makebox(0,0){$c_{1}$}}
\put(2.5,6){\makebox(0,0){$c$}} \put(2.25,-.25){\line(0,1){6.5}}
\put(9.7,-.25){\line(0,1){6.5}} \put(2.8,-.25){\line(0,1){6.5}}
\put(6.6,-.25){\line(0,1){6.5}}
\put(5.05,-.25){\line(0,1){6.5}} \put(8.05,-.25){\line(0,1){6.5}}
\put(1.4,5.7){\line(1,0){8.3}} \put(1.4,-.25){\line(1,0){8.3}}
\put(1.4,.95){\line(1,0){8.3}}
\put(1.4,2.15){\line(1,0){8.3}}\put(1.4,3.35){\line(1,0){8.3}}
\put(1.4,5.15){\line(1,0){8.3}}
\put(3.25,5.4){\makebox(0,0){$
1$}} \put(4,5.4){\makebox(0,0){$
1$}}\put(4.75,5.4){\makebox(0,0){$ 1$}}

\put(2.5,4.8){\makebox(0,0){$ 1$}} \put(2.5,4.2){\makebox(0,0){$
1$}}\put(2.5,3.6){\makebox(0,0){$ 1$}}
\put(5.5,4.8){\makebox(0,0){$ 1$}}
\put(6.25,4.8){\makebox(0,0){$ 1$}}
\put(7,4.2){\makebox(0,0){$ 1$}}
\put(7.75,4.2){\makebox(0,0){$ 1$}}
\put(8.5,3.6){\makebox(0,0){$ 1$}}
\put(9.25,3.6){\makebox(0,0){$ 1$}}
\put(3.25,3){\makebox(0,0){$ 1$}}
\put(3.25,2.4){\makebox(0,0){$ 1$}}
\put(4,1.8){\makebox(0,0){$ 1$}}
\put(4,1.2){\makebox(0,0){$ 1$}}
\put(4.75,.6){\makebox(0,0){$ 1$}}
\put(4.75,0){\makebox(0,0){$ 1$}}

\put(6.25,3){\makebox(0,0){$ 1$}}
\put(9.25,3){\makebox(0,0){$ 1$}}
\put(5.5,2.4){\makebox(0,0){$ 1$}}
\put(7,2.4){\makebox(0,0){$ 1$}}
\put(6.25,1.8){\makebox(0,0){$ 1$}}
\put(7.75,1.8){\makebox(0,0){$ 1$}}

\put(7,1.2){\makebox(0,0){$ 1$}}
\put(8.5,1.2){\makebox(0,0){$ 1$}}
\put(7.75,.6){\makebox(0,0){$ 1$}}
\put(9.25,.6){\makebox(0,0){$ 1$}}
\put(5.5,0){\makebox(0,0){$ 1$}}
\put(8.5,0){\makebox(0,0){$ 1$}}
\put(.5,3){\makebox(0,0){$A_2\,${\rm :}}}

\put(15,5.5){\circle*{.3}} \put(17,5.5){\circle*{.3}}
\put(13.7,3.5){\circle*{.3}}\put(18.3,3.5){\circle*{.3}}
\put(15,1.5){\circle*{.3}}\put(17,1.5){\circle*{.3}}
\put(16,3.5){\circle*{.3}} \put(15,4.3){\circle*{.3}}
\put(17,4.3){\circle*{.3}} \put(16,2.5){\circle*{.3}}

\put(15,5.5){\line(1,0){2}} \put(15,5.5){\line(-2,-3){1.4}}
\put(15,5.5){\line(0,-1){1.2}} \put(15,1.5){\line(1,0){2}}
\put(15,1.5){\line(-2,3){1.4}} \put(15,1.5){\line(1,1){1}}
\put(17,5.5){\line(0,-1){1.2}} \put(17,5.5){\line(2,-3){1.4}}
\put(17,1.5){\line(-1,1){1}}\put(17,1.5){\line(2,3){1.4}}
\put(16,3.5){\line(0,-1){1}} \put(16,3.5){\line(-6,5){1}}
\put(16,3.5){\line(6,5){1}}
\put(13.7,3.5){\line(3,2){1.3}}\put(18.3,3.5){\line(-3,2){1.3}}

\put(16.4,3.5){\makebox(0,0){$c$}}
\put(17.8,1.5){\makebox(0,0){$c_{21}$}}
\put(13,3.5){\makebox(0,0){$c_{31}$}}
\put(17.6,5.7){\makebox(0,0){$c_{11}$}}
\put(19,3.5){\makebox(0,0){$c_{12}$}}
\put(14.4,5.7){\makebox(0,0){$c_{32}$}}
\put(14.3,1.5){\makebox(0,0){$c_{22}$}}
\put(16.5,2.7){\makebox(0,0){$c_2$}}
\put(17.2,3.8){\makebox(0,0){$c_1$}}
\put(15,3.8){\makebox(0,0){$c_3$}}

\put(16.3,0){\makebox(0,0){$T_2$}}

\end{picture}\medskip

 \begin{theorem}\label{resultforkeq3}
Let  $A \in {\frak D}_{10,3}$. Then the following are equivalent:
$$
\begin{array}{rl}
(i) & \Theta^3(A) = A\,, \;\makebox{but}\; \Theta(A)
\ne A\,, \\&\;\makebox{i.e.} \;A\;\makebox{is a proper
solution of} \;(gHS)\;\makebox{with} \; m = 3\,;\\
(ii) & A\; \makebox{is an incidence matrix of the} \;
10_3 \; \makebox{configuration} \; 10_3F\,;\\
(iii) & A\; \makebox{is an adjacency matrix of
the Terwilliger graph} \; T_2\,;\\
 (iv) & A \; \makebox{is p--equivalent to the
 HS--form}\; A_2\,.
\end{array}
$$
\end{theorem}

\Prf This is an immediate consequence of \cite[Proposition
4.2]{AFLN1}.\qed

 \begin{theorem}\label{resultforkeq4}
Let  $A \in {\frak D}_{17,4}$. Then the following are equivalent:
$$
\begin{array}{rl}
(i) & \Theta^2(A) = A\,, \;\makebox{i.e.} \;A\;\makebox{is a solution of} \;(gHS)\;\makebox{with} \; m = 2\,;\\
(ii) & A\; \makebox{is an incidence matrix of the} \; 17_4 \; \makebox{configuration} \;\# 1971 \, \makebox{in Betten's list \cite{BB99}} ;\\
(iii) & A \; \makebox{is p--equivalent to the HS--form}\;
S_{HS}(P)\,, where
\end{array}
$$
$$
P \; = \;\left(\begin{array}{ccc|ccc|ccc|ccc}
&&&1&&&1&&&1&&\\
&&&&1&&&1&&&1&\\
&&&&&1&&&1&&&1\\\hline
1&&&&&&1&&&&1&\\
&1&&&&&&1&&&&1\\
&&1&&&&&&1&1&&\\\hline
1&&&1&&&&&&&&1\\
&1&&&1&&&&&1&&\\
&&1&&&1&&&&&1&\\\hline
1&&&&&1&&1&&&&\\
&1&&1&&&&&1&&&\\
&&1&&1&&1&&&&&\\
\end{array}\right)
$$
\end{theorem}

\Prf In Proposition \ref{174} we showed that \cite[Conjecture
4.3]{AFLN1} is in fact a theorem. Then the statements follow by
applying \cite[Theorem 4.7]{AFLN1}. \qed

\begin{remark}
The solutions of $(gHS)$ for $m=1$ yield graphs which can be seen as association schemes.
However, this is not the case in general. In fact, for $m=3$, the graph $T_2$
from Theorem \ref{resultforkeq3}(iii), cannot be seen as an association scheme
since its vertex $c$ is the only one not lying in a $3$--cycle.
\end{remark}

\end{document}